\documentclass[12pt]{article}
\usepackage{latexsym,amssymb,amsmath,amsthm,enumerate,geometry,cite,float,tikz}
\newtheorem{theorem}{Theorem}

\usepackage{lineno}
\usepackage{setspace}
\allowdisplaybreaks

\begin{document}
\onehalfspace

\title{Relating dissociation, independence, and matchings}
\author{Felix Bock\and Johannes Pardey\and Lucia D. Penso\and Dieter Rautenbach}
\date{}

\maketitle
\vspace{-10mm}
\begin{center}
{\small 
Institute of Optimization and Operations Research, Ulm University, Ulm, Germany\\
\texttt{$\{$felix.bock,johannes.pardey,lucia.penso,dieter.rautenbach$\}$@uni-ulm.de}
}
\end{center}

\begin{abstract}
A dissociation set in a graph is a set of vertices 
inducing a subgraph of maximum degree at most $1$. 
Computing the dissociation number ${\rm diss}(G)$ of a given graph $G$,
defined as the order of a maximum dissociation set in $G$,
is algorithmically hard even when $G$ is restricted to be bipartite.
Recently, Hosseinian and Butenko 
proposed a simple $\frac{4}{3}$-approximation algorithm 
for the dissociation number problem in bipartite graphs.
Their result relies on the inequality
${\rm diss}(G)\leq\frac{4}{3}\alpha(G-M)$
implicit in their work,
where $G$ is a bipartite graph,
$M$ is a maximum matching in $G$,
and $\alpha(G-M)$ denotes the independence number of $G-M$.
We show that the pairs $(G,M)$ 
for which this inequality holds with equality
can be recognized efficiently,
and that a maximum dissociation set can be determined for them efficiently.
The dissociation number of a graph $G$ satisfies 
$\max\{ \alpha(G),2\nu_s(G)\}
\leq {\rm diss}(G)\leq \alpha(G)+\nu_s(G)\leq 2\alpha(G)$,
where $\nu_s(G)$ denotes the induced matching number of $G$.
We show that deciding whether ${\rm diss}(G)$ 
equals any of the four terms 
lower and upper bounding ${\rm diss}(G)$ is NP-hard.\\[3mm]
{\bf Keywords:} Dissociation set; independent set; matching; induced matching
\end{abstract}

\pagebreak

\section{Introduction}
We consider finite, simple, and undirected graphs, and use standard terminology.
A set $I$ of vertices of a graph $G$
is a {\it dissociation set} in $G$ if the subgraph $G[I]$ of $G$ induced by $I$
has maximum degree at most $1$, 
and the {\it dissociation number ${\rm diss}(G)$} of $G$ 
is the order of a maximum dissociation set in $G$.
Dissociation sets and the dissociation number were introduced as a special vertex-deletion problem by Yannakakis \cite{ya} 
who showed that the dissociation number problem,
that is, the problem of deciding whether ${\rm diss}(G)\geq k$
for a given pair $(G,k)$, where $G$ is a graph and $k$ is a positive integer,
is NP-complete even restricted to instances where $G$ is a bipartite graph.
This initial hardness result was stengthened in different ways \cite{bocalo,ordofigowe};
in particular, 
the problem remains NP-complete for bipartite graphs of maximum degree at most $3$.

Recently, Hosseinian and Butenko \cite{hobu} 
proposed a simple $\frac{4}{3}$-approximation algorithm 
for the maximum dissociation set problem restricted to bipartite graphs.
Their result can be derived from the following two simple inequalities:
Let $G$ be a graph and let $M$ be a maximum matching in $G$.
Since every independent set in $G-M$ is a dissociation set in $G$, we have
\begin{eqnarray}
{\rm diss}(G) & \geq & \alpha(G-M)\label{e1},
\end{eqnarray}
where $\alpha(H)$ denotes the {\it independence number} of a graph $H$,
which is the order of a maximum independent set in $H$.
Now, if $G$ is bipartite, then one can show
\begin{eqnarray}
{\rm diss}(G) & \leq & \frac{4}{3}\alpha(G-M)\label{e2}.
\end{eqnarray}
Since a maximum matching $M$ in a given bipartite graph $G$
as well as a maximum independent set $I$ in the bipartite graph $G-M$
can be determined efficiently, the combination of (\ref{e1}) and (\ref{e2}) 
implies that returning $I$ yields a $\frac{4}{3}$-approximation 
for the maximum dissociation set problem in $G$.
We will give the simple proof of (\ref{e2}) that is implicit in \cite{hobu} below.
As our first contribution we show that the extremal graphs for (\ref{e2})
have a very restricted structure, which yields the following.

\begin{theorem}\label{theorem1}
For a given pair $(G,M)$, 
where $G$ is a bipartite graph and $M$ is a maximum matching in $G$, 
one can decide in polynomial time whether (\ref{e2}) is satisfied with equality. Furthermore, in case of equality, 
one can determine in polynomial time a maximum dissociation set in $G$.
\end{theorem}
Next to (\ref{e1}) and (\ref{e2}), 
there are the following relations between 
the dissociation number ${\rm diss}(G)$,
the independence number $\alpha(G)$,
and the induced matching number $\nu_s(G)$ of a graph $G$:
\begin{eqnarray}
{\rm diss}(G) & \leq &2\alpha(G),\label{e3}\\
{\rm diss}(G) & \geq &2\nu_s(G),\label{e3b}\\
{\rm diss}(G) & \geq &\alpha(G),\mbox{ and}\label{e3c}\\
{\rm diss}(G) & \leq &\alpha(G)+\nu_s(G).\label{e3d}
\end{eqnarray}
While these inequalities are all straightforward, 
the extremal graphs are not easy to describe,
and we show the following.

\begin{theorem}\label{theorem2}
For each of the inequalities (\ref{e3}), (\ref{e3b}), (\ref{e3c}), and (\ref{e3d}), 
it is NP-hard to decide whether a given graph satisfies it with equality.
\end{theorem}
In view of the special role of bipartite graphs in this context,
it makes sense to consider the bipartite extremal graphs
for (\ref{e3}) to (\ref{e3d}).
It is easy to see that a bipartite graph $G$ satisfies 
${\rm diss}(G)=2\alpha(G)$
if and only if $G$ is $1$-regular.
For a bipartite graph $G$,
the equality ${\rm diss}(G)=\alpha(G)$ holds if and only if
$G$ has no induced matching $M$ intersecting 
every maximum matching in $G$.
Inspecting the proofs in \cite{zeripiwecobe} reveals 
that Zenklusen et al.~showed that it its NP-complete to decide, for a given pair $(G,k)$, 
where $G$ is a bipartite graph and $k$ is a positive integer,
whether there is an induced matching $M$ in $G$ of size $|M|$ at most $k$ 
intersecting every maximum matching in $G$.
Unfortunately, the size bound is crucial for their reduction.
The complexity of the induced matching number
is closely tied to the complexity of the dissociation number
\cite{ordofigowe}. 
In particular, it is hard for bipartite graphs,
and the complexity of recognizing the bipartite extremal graphs for (\ref{e3b}) and (\ref{e3d}) is open.
The close relation between dissociation sets, independent sets,
and (induced) matchings also reflects in the obvious relation 
$${\rm diss}(G)=\max\{ \alpha(G-M):M\mbox{ is an induced matching in }G\}.$$
Before we proceed to the proofs of our results, 
we briefly mention related research.
In fact, bounds on the dissociation number \cite{brkakase,brjakaseta},
fast exact algorithms \cite{kakasc,xiko}, 
randomized approximation algorithms \cite{kakasc}, and 
fixed parameter tractability \cite{ts}
were studied.
As observed in several references,
dissociations sets are the dual of so-called 
{\it $3$-path (vertex) covers}, cf.~also \cite{beocra}.

\section{Proofs}

The following two subsections contain the proofs of Theorem \ref{theorem1} and Theorem \ref{theorem2}.

\subsection{Structure and recognition of the extremal graphs for (\ref{e2})}

We first give a proof of (\ref{e2}), which is implicit in \cite{hobu}.
After that we consider the extremal graphs.
Throughout this subsection, let $G$ be a bipartite graph of order $n$ 
with partite sets $A$ and $B$,
and let $M$ be a maximum matching in $G$.
Let $I$ be a maximum dissociation set in $G$.
Let $E$ be the induced matching spanned by $I$, that is, $E=E(G[I])$.
Note that $|E|\leq \frac{|I|}{2}$.
Let $A=A_1\cup A_2\cup A_3\cup A_4$
and $B=B_1\cup B_2\cup B_3\cup B_4$ be partitions of $A$ and $B$ such that 
\begin{itemize}
\item $E\cap M$ is a perfect matching between $A_1$ and $B_1$,
\item $E\setminus M$ is a perfect matching between $A_2$ and $B_2$, 
\item $A_3\cup B_3$ is the set of isolated vertices in $G[I]$, 
\item $A_4=A\setminus I$, and $B_4=B\setminus I$.
\end{itemize}
Let $\ell=|A_2|$.
Since $I\setminus A_2$ is an independent set in $G-M$, we have
\begin{eqnarray}\label{e5}
\alpha(G-M) &\geq & |I\setminus A_2|=|I|-|A_2|={\rm diss}(G)-\ell.
\end{eqnarray}
Gallai's theorem implies that 
$\alpha(G)$ and the vertex cover number $\tau(G)$ of $G$ add up to the order $n$ of $G$, and
K\H{o}nig's theorem implies that 
$\tau(G)$ equals the matching number $|M|$ of $G$, 
which together implies that $\alpha(G)=n-|M|$.
Since $M$ contains $|E\cap M|=|E|-\ell\leq \frac{|I|}{2}-\ell$ 
edges spanned by $I$,
and at most $n-|I|$ further edges, 
one incident with each vertex of $A_4\cup B_4=V(G)\setminus I$, we have 
\begin{eqnarray}\label{e6}
|M| &\leq & |E\cap M|+(n-|I|)
\leq \left(\frac{|I|}{2}-\ell\right)+\left(n-|I|\right)
=n-\ell-\frac{{\rm diss}(G)}{2},
\end{eqnarray}
and, hence,
\begin{eqnarray}\label{e4}
\alpha(G-M)
\geq \alpha(G)
=n-|M|
\stackrel{(\ref{e6})}{\geq} n-\left(n-\ell-\frac{{\rm diss}(G)}{2}\right)
=\ell+\frac{{\rm diss}(G)}{2}.
\end{eqnarray}
Adding (\ref{e5}) and (\ref{e4}) yields (\ref{e2}).

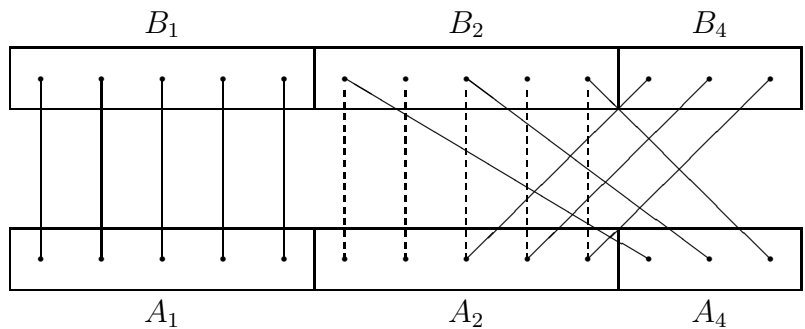
\begin{figure}[h]
\begin{center}
\unitlength 0.8mm 
\linethickness{0.4pt}
\ifx\plotpoint\undefined\newsavebox{\plotpoint}\fi 
\begin{picture}(135,49)(0,0)
\put(130,10){\circle*{1}}
\put(130,40){\circle*{1}}
\put(120,10){\circle*{1}}
\put(120,40){\circle*{1}}
\put(110,10){\circle*{1}}
\put(110,40){\circle*{1}}
\put(100,10){\circle*{1}}
\put(100,40){\circle*{1}}
\put(90,10){\circle*{1}}
\put(90,40){\circle*{1}}
\put(80,10){\circle*{1}}
\put(80,40){\circle*{1}}
\put(70,10){\circle*{1}}
\put(70,40){\circle*{1}}
\put(60,10){\circle*{1}}
\put(60,40){\circle*{1}}
\put(50,10){\circle*{1}}
\put(50,40){\circle*{1}}
\put(40,10){\circle*{1}}
\put(40,40){\circle*{1}}
\put(30,10){\circle*{1}}
\put(30,40){\circle*{1}}
\put(20,10){\circle*{1}}
\put(20,40){\circle*{1}}
\put(10,10){\circle*{1}}
\put(10,40){\circle*{1}}
\put(105,5){\framebox(30,10)[]{}}
\put(105,35){\framebox(30,10)[]{}}
\put(55,5){\framebox(50,10)[]{}}
\put(55,35){\framebox(50,10)[]{}}
\put(5,5){\framebox(50,10)[]{}}
\put(5,35){\framebox(50,10)[]{}}
\put(120,1){\makebox(0,0)[]{$A_4$}}
\put(120,49){\makebox(0,0)[]{$B_4$}}
\put(80,1){\makebox(0,0)[]{$A_2$}}
\put(30,1){\makebox(0,0)[]{$A_1$}}
\put(80,49){\makebox(0,0)[]{$B_2$}}
\put(30,49){\makebox(0,0)[]{$B_1$}}
\put(50,10){\line(0,1){30}}
\put(40,10){\line(0,1){30}}
\put(30,10){\line(0,1){30}}
\put(20,10){\line(0,1){30}}
\put(10,10){\line(0,1){30}}
\put(99.93,9.93){\line(0,1){.9667}}
\put(99.93,11.863){\line(0,1){.9667}}
\put(99.93,13.796){\line(0,1){.9667}}
\put(99.93,15.73){\line(0,1){.9667}}
\put(99.93,17.663){\line(0,1){.9667}}
\put(99.93,19.596){\line(0,1){.9667}}
\put(99.93,21.53){\line(0,1){.9667}}
\put(99.93,23.463){\line(0,1){.9667}}
\put(99.93,25.396){\line(0,1){.9667}}
\put(99.93,27.33){\line(0,1){.9667}}
\put(99.93,29.263){\line(0,1){.9667}}
\put(99.93,31.196){\line(0,1){.9667}}
\put(99.93,33.13){\line(0,1){.9667}}
\put(99.93,35.063){\line(0,1){.9667}}
\put(99.93,36.996){\line(0,1){.9667}}
\put(89.93,9.93){\line(0,1){.9667}}
\put(89.93,11.863){\line(0,1){.9667}}
\put(89.93,13.796){\line(0,1){.9667}}
\put(89.93,15.73){\line(0,1){.9667}}
\put(89.93,17.663){\line(0,1){.9667}}
\put(89.93,19.596){\line(0,1){.9667}}
\put(89.93,21.53){\line(0,1){.9667}}
\put(89.93,23.463){\line(0,1){.9667}}
\put(89.93,25.396){\line(0,1){.9667}}
\put(89.93,27.33){\line(0,1){.9667}}
\put(89.93,29.263){\line(0,1){.9667}}
\put(89.93,31.196){\line(0,1){.9667}}
\put(89.93,33.13){\line(0,1){.9667}}
\put(89.93,35.063){\line(0,1){.9667}}
\put(89.93,36.996){\line(0,1){.9667}}
\put(79.93,9.93){\line(0,1){.9667}}
\put(79.93,11.863){\line(0,1){.9667}}
\put(79.93,13.796){\line(0,1){.9667}}
\put(79.93,15.73){\line(0,1){.9667}}
\put(79.93,17.663){\line(0,1){.9667}}
\put(79.93,19.596){\line(0,1){.9667}}
\put(79.93,21.53){\line(0,1){.9667}}
\put(79.93,23.463){\line(0,1){.9667}}
\put(79.93,25.396){\line(0,1){.9667}}
\put(79.93,27.33){\line(0,1){.9667}}
\put(79.93,29.263){\line(0,1){.9667}}
\put(79.93,31.196){\line(0,1){.9667}}
\put(79.93,33.13){\line(0,1){.9667}}
\put(79.93,35.063){\line(0,1){.9667}}
\put(79.93,36.996){\line(0,1){.9667}}
\put(69.93,9.93){\line(0,1){.9667}}
\put(69.93,11.863){\line(0,1){.9667}}
\put(69.93,13.796){\line(0,1){.9667}}
\put(69.93,15.73){\line(0,1){.9667}}
\put(69.93,17.663){\line(0,1){.9667}}
\put(69.93,19.596){\line(0,1){.9667}}
\put(69.93,21.53){\line(0,1){.9667}}
\put(69.93,23.463){\line(0,1){.9667}}
\put(69.93,25.396){\line(0,1){.9667}}
\put(69.93,27.33){\line(0,1){.9667}}
\put(69.93,29.263){\line(0,1){.9667}}
\put(69.93,31.196){\line(0,1){.9667}}
\put(69.93,33.13){\line(0,1){.9667}}
\put(69.93,35.063){\line(0,1){.9667}}
\put(69.93,36.996){\line(0,1){.9667}}
\put(59.93,9.93){\line(0,1){.9667}}
\put(59.93,11.863){\line(0,1){.9667}}
\put(59.93,13.796){\line(0,1){.9667}}
\put(59.93,15.73){\line(0,1){.9667}}
\put(59.93,17.663){\line(0,1){.9667}}
\put(59.93,19.596){\line(0,1){.9667}}
\put(59.93,21.53){\line(0,1){.9667}}
\put(59.93,23.463){\line(0,1){.9667}}
\put(59.93,25.396){\line(0,1){.9667}}
\put(59.93,27.33){\line(0,1){.9667}}
\put(59.93,29.263){\line(0,1){.9667}}
\put(59.93,31.196){\line(0,1){.9667}}
\put(59.93,33.13){\line(0,1){.9667}}
\put(59.93,35.063){\line(0,1){.9667}}
\put(59.93,36.996){\line(0,1){.9667}}
\put(130,40){\line(-1,-1){30}}
\put(130,10){\line(-1,1){30}}
\put(120,10){\line(-4,3){40}}
\put(120,40){\line(-1,-1){30}}
\put(110,10){\line(-5,3){50}}
\put(110,40){\line(-1,-1){30}}
\end{picture}
\caption{The continous lines illustrate the edges in $M$ 
while the dashed lines illustrate those in $E\setminus M$.
All remaining edges of $G$ are not illustrated 
and intersect $A_4\cup B_4$.}\label{fig2}
\end{center}
\end{figure}

We now consider the extremal graphs for (\ref{e2}),
which leads to a proof of Theorem \ref{theorem1}.
Therefore, we suppose that (\ref{e2}) is satisfied with equality.
This implies that equality holds throughout the inequality chains 
(\ref{e5}), (\ref{e6}), and (\ref{e4}).
Equality throughout (\ref{e5}) and (\ref{e4}) implies
$\ell=\frac{{\rm diss}(G)}{4}$ and $\alpha(G-M)=\alpha(G)$.
Equality in (\ref{e6}) implies 
$|E|
=|E\cap M|+|E\setminus M|
=|E\cap M|+\ell
=\frac{|I|}{2}$,
which implies that $G[I]$ is $1$-regular, or, equivalently, $A_3=B_3=\emptyset$.
In view of the above value of $\ell$,
exactly half the edges of $G[I]$ belong to $M$, or, equivalently,
$\ell=|A_1|=|A_2|=|B_1|=|B_2|$.
Furthermore, equality in (\ref{e6}) also implies that the matching $M$ 
contains exactly $n-|I|=|A_4|+|B_4|$ further edges,
one incident with each vertex of $A_4\cup B_4$;
these edges match all of $A_4$ into $B_2$ as well as all of $B_4$ into $A_2$, 
in particular, $\ell\geq |A_4|,|B_4|$.
The matching $M$ leaves exactly 
$(|A_2|-|B_4|)+(|B_2|-|A_4|)=2\ell-(|A_4|+|B_4|)$ 
vertices unmatched that all lie in $A_2\cup B_2$.
See Figure \ref{fig2} for an illustration.

\begin{figure}[h]
\begin{center}
\unitlength 1mm 
\linethickness{0.4pt}
\ifx\plotpoint\undefined\newsavebox{\plotpoint}\fi 
\begin{picture}(48,48)(0,0)
\put(15,5){\circle*{1}}
\put(15,45){\circle*{1}}
\put(35,5){\circle*{1}}
\put(35,45){\circle*{1}}
\put(45,25){\circle*{1}}
\put(5,25){\circle*{1}}
\put(15,45){\line(1,0){20}}
\put(35,45){\line(1,-2){10}}
\put(45,25){\line(-1,-2){10}}
\put(35,5){\line(-1,0){20}}
\put(15,5){\line(-1,2){10}}
\put(5,25){\line(1,2){10}}
\put(15,48){\makebox(0,0)[cc]{$A_1$}}
\put(35,48){\makebox(0,0)[cc]{$B_1$}}
\put(48,25){\makebox(0,0)[cc]{$A_4$}}
\put(2,25){\makebox(0,0)[cc]{$B_4$}}
\put(35,2){\makebox(0,0)[cc]{$B_2$}}
\put(15,2){\makebox(0,0)[cc]{$A_2$}}
\put(25,48){\makebox(0,0)[cc]{$M$}}
\put(43,14){\makebox(0,0)[cc]{$M$}}
\put(7,14){\makebox(0,0)[cc]{$M$}}
\put(43,36){\makebox(0,0)[cc]{$M'$}}
\put(7,36){\makebox(0,0)[cc]{$M'$}}
\put(25,2){\makebox(0,0)[cc]{$M'$}}
\end{picture}
\caption{The cyclic order respected by the components of $H$,
that is, a cycle in $H$ traverses the sets in the cyclic order
$A_1,B_1,A_4,B_2,A_2,B_4,A_1,B_1,A_4,B_2,A_2,B_4,\ldots$.}\label{fig1}
\end{center}
\end{figure}
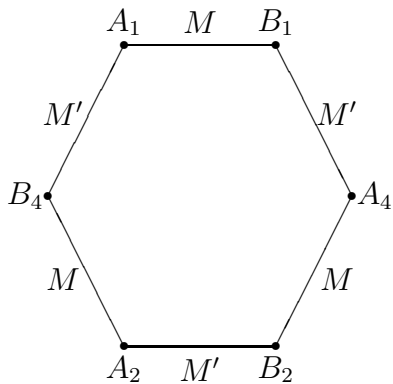

Now, let $M'$ be any maximum matching in $G-M$.
Since $\alpha(G-M)=\alpha(G)$, the results of Gallai and K\H{o}nig 
imply that $G-M$ and $G$ have the same matching number, 
which implies $|M|=|M'|$.
Since in $G-M$ the vertices in $A_1$ have all their neighbors in $B_4$,
and the vertices in $B_1$ have all their neighbors in $A_4$,
the matching $M'$ leaves at least 
$(|A_1|-|B_4|)+(|B_1|-|A_4|)=2\ell-(|A_4|+|B_4|)$ vertices 
in $A_1\cup B_1$ unmatched.
Since $|M|=|M'|$, 
this actually implies that $M'$ consists of 
$|B_4|$ edges matching all of $B_4$ into $A_1$,
$|A_4|$ edges matching all of $A_4$ into $B_1$, and 
the $\ell$ edges from $E\setminus M$ that form a perfect matching between $A_2$ and $B_2$.
Let $H$ be the graph with vertex set $V(G)$ and edge set $M\cup M'$.
The components of $H$ are $M$-$M'$-alternating paths and cycles
that traverse the sets $A_i$ and $B_i$ 
respecting the cyclic order illustrated in Figure \ref{fig1}.

For components $C$ of $H$ that are cycles, 
this implies that the length of $C$ is a multiple of $6$.
Furthermore, 
there are exactly $2\ell-(|A_4|+|B_4|)$ ccomponents of $H$ that are paths;
they all have length $4$ modulo $6$, and,
starting with an edge from $M$,
they go from $A_1\cup B_1$ to $A_2\cup B_2$.
For a component $P$ of $H$ that is a path, these structural properties 
allow to decide the location of the individual vertices.
If, for example, the path $P:u_1u_2\ldots u_{12}\ldots$ starts in a vertex $u_1$ from $A$ not covered by $M'$,
then 
$u_1,u_7,\ldots  \in A_1$,
$u_2,u_8,\ldots  \in B_1$,
$u_3,u_9,\ldots  \in A_4$,
$u_4,u_{10},\ldots  \in B_2$,
$u_5,u_{11},\ldots  \in A_2$, and
$u_6,u_{12},\ldots  \in B_4$.

We now formulate a satisfiable $2$-{\sc Sat} formula $f$
such that a satisfying truth assignment allows to derive a (suitable) 
location of all vertices of $G$ on cycle components of $H$.
Therefore, let 
$C_1,\ldots,C_k$ be the components of $H$ that are cycles.
For $i$ in $[k]$,
let $C_i:a_i^1b_i^1a_i^2b_i^2\ldots b_i^{n_i}a_i^1$,
where $a_i^1\in A$ and $a_i^1b_i^1\in M$.
Note that exactly one of the three vertices $a_i^1$, $a_i^2$, and $a_i^3$ belongs to $A_4$, which also determines the location of every other vertex on $C_i$.
For every $i$ in $[k]$, 
we introduce three Boolean variables $x_i^1$, $x_i^2$, and $x_i^3$, and 
add to $f$ the three clauses
$\bar{x}_i^1\vee \bar{x}_i^2$,
$\bar{x}_i^2\vee \bar{x}_i^3$, and
$\bar{x}_i^1\vee \bar{x}_i^3$;
where $x_i^j$ being true corresponds to $a_i^j\in A_4$.
Now, we consider an edge $ab$ of $G$ that does not belong to $M\cup M'$
such that at least one endpoint of $ab$ lies on a cycle component of $H$.
The structural properties imply that $a\in A_4$ or $b\in B_4$ (or both).
If $a$ lies on a cycle component $C$ of $H$, 
$b$ lies on a path component of $H$, and 
$b\not\in B_4$, 
then we add to $f$ the clause $x_i^j$, 
where $i$ and $j$ are such that $C$ is $C_i$, and 
the distance of $a_i^j$ and $a$ on $C$ is $0$ modulo $6$.
Note that, if $b\in B_4$, then no clause is added to $f$.
Similarly, 
if $a$ lies on a path component of $H$, 
$b$ lies on a cycle component of $H$, 
and $a\not\in A_4$, then we proceed analogously 
by adding to $f$ the clause $x_i^j$
corresponding to the condition $x_i^j\in A_4$ 
that is equivalent to the condition $b\in B_4$.
Finally, if $a$ and $b$ both lie on cycle components of $H$,
say $a=a_i^j$ and $b=b_{i'}^{j'}$,
then we add to $f$ the clause $x_i^j\vee x_{i'}^{j''}$,
where $j''$ is the uniquely determined index such that 
$a_{i'}^{j''}$ is the unique vertex in $\{ a_{i'}^1,a_{i'}^2,a_{i'}^3\}$ 
whose distance to $b_{i'}^{j'}\in B_4$ within $C_{i'}$
is equivalent to $3$ modulo $6$.
This completes the construction of $f$.
Clearly, setting $x_i^j$ to true for every $i$ and $j$ with $a_i^j\in A_4$
and false otherwise, yields a satisfying truth assignment for $f$,
that is, if (\ref{e2}) is satisfied with equality, then $f$ is satisfiable.

Given a bipartite graph $G$ and a maximum matching $M$ in $G$,
one can in polynomial time
\begin{itemize}
\item determine a bipartition $A$ and $B$ of $G$,
\item determine a maximum matching $M'$ of $G-M$,
\item check the necessary condition $|M|=|M'|$,
\item construct the auxiliary graph $H$, 
\item check the necessary condition that all cycle components of $H$
have length $0$ modulo $6$,
and all path components of $H$ have length $4$ modulo $6$,
\item suitably assign the vertices in path components of $H$ 
to the sets $A_i$ and $B_i$ with $i\in \{ 1,2,4\}$ as described above,
\item check the necessary condition that all edges of $G$ 
that do not belong to $M\cup M'$
and connect vertices in path components of $H$
intersect $A_4\cup B_4$, and 
\item set up the $2$-{\sc Sat} formula $f$ as described above,
check its satisfiability, and, in case of satisfiability,
determine a satisfying truth assignment.
\end{itemize}
Note that the construction of $f$ requires only available knowledge.
As we have seen above, 
if (\ref{e2}) is satisfied with equality, 
then the necessary conditions mentioned above hold, 
and $f$ is satisfiable.
Conversely, 
if the necessary conditions mentioned above hold, and $f$ is satisfiable,
then a satisfying truth assignment allows to suitably assign 
the vertices in cycle components of $H$ 
to the sets $A_i$ and $B_i$ with $i\in \{ 1,2,4\}$
in such a way that 
all edges of $G$ that do not belong to $M\cup M'$ intersect $A_4\cup B_4$,
$A_1\cup A_2\cup B_1$ is a maximum independent set in $G-M$, and 
$A_1\cup A_2\cup B_1\cup B_2$ is a dissociation set in $G$,
that is, (\ref{e2}) is satisfied with equality.
Note that $H$ may contain cycle components $C_i$
for which all three variables $x_i^1$, $x_i^2$, and $x_i^3$ are false, 
even if $G$ contains edges connecting $C_i$ 
to other components of $H$.
In such a case, any of the three vertices $a_i^1$, $a_i^2$, and $a_i^3$
may be located within $A_4$, which yields three different valid possibilities.
Note furthermore that, in case of equality in (\ref{e2}),
the set $A_1\cup A_2\cup B_1\cup B_2$, 
which is efficiently constructible as explained above,
is a maximum dissociation set in $G$.

This complete the proof of Theorem \ref{theorem1}.

\subsection{Hardness of deciding equality in (\ref{e3}) to (\ref{e3d})}

In this subsection, 
we show Theorem \ref{theorem2}.

For the hardness of deciding equality in (\ref{e3}), (\ref{e3b}), or (\ref{e3c}),
we suitably adapt Karp's proof \cite{ka} 
of the NP-completeness of the {\sc Clique} problem,
reducing $3$-{\sc Sat} to the respective problems.
Therefore, 
let $f$ be an instance of $3$-{\sc Sat}
consisting of the clauses $C_1,\ldots,C_m$
over the Boolean variables $x_1,\ldots,x_n$.

\medskip

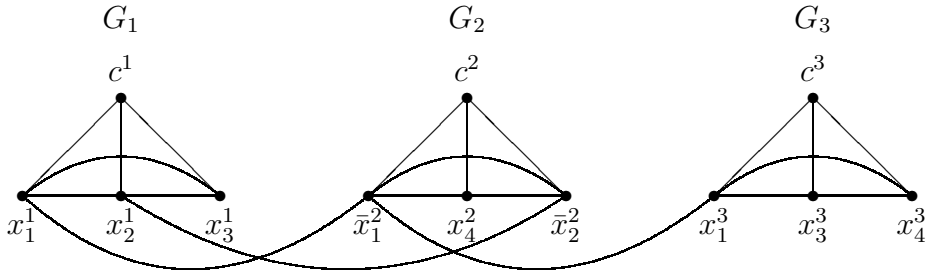
\begin{figure}[h]
\begin{center}
\unitlength 1.3mm 
\linethickness{0.4pt}
\ifx\plotpoint\undefined\newsavebox{\plotpoint}\fi 
\begin{picture}(96,28)(0,0)
\put(5,10){\circle*{1}}
\put(40,10){\circle*{1}}
\put(75,10){\circle*{1}}
\put(15,10){\circle*{1}}
\put(50,10){\circle*{1}}
\put(85,10){\circle*{1}}
\put(25,10){\circle*{1}}
\put(60,10){\circle*{1}}
\put(95,10){\circle*{1}}
\put(15,20){\circle*{1}}
\put(50,20){\circle*{1}}
\put(85,20){\circle*{1}}
\put(15,20){\line(-1,-1){10}}
\put(50,20){\line(-1,-1){10}}
\put(85,20){\line(-1,-1){10}}
\put(5,10){\line(1,0){20}}
\put(40,10){\line(1,0){20}}
\put(75,10){\line(1,0){20}}
\put(25,10){\line(-1,1){10}}
\put(60,10){\line(-1,1){10}}
\put(95,10){\line(-1,1){10}}
\put(15,20){\line(0,-1){10}}
\put(50,20){\line(0,-1){10}}
\put(85,20){\line(0,-1){10}}
\qbezier(25,10)(15,18)(5,10)
\qbezier(60,10)(50,18)(40,10)
\qbezier(95,10)(85,18)(75,10)
\put(15,28){\makebox(0,0)[cc]{$G_1$}}
\put(50,28){\makebox(0,0)[cc]{$G_2$}}
\put(85,28){\makebox(0,0)[cc]{$G_3$}}
\put(15,23){\makebox(0,0)[cc]{$c^1$}}
\put(50,23){\makebox(0,0)[cc]{$c^2$}}
\put(85,23){\makebox(0,0)[cc]{$c^3$}}
\qbezier(5,10)(21.5,-5)(40,10)
\qbezier(40,10)(56.5,-5)(75,10)
\put(5,7){\makebox(0,0)[cc]{$x^1_1$}}
\put(40,7){\makebox(0,0)[cc]{$\bar{x}^2_1$}}
\put(75,7){\makebox(0,0)[cc]{$x^3_1$}}
\put(15,7){\makebox(0,0)[cc]{$x^1_2$}}
\put(60,7){\makebox(0,0)[cc]{$\bar{x}^2_2$}}
\put(49.5,7){\makebox(0,0)[cc]{$x^2_4$}}
\put(85,7){\makebox(0,0)[cc]{$x^3_3$}}
\put(95,7){\makebox(0,0)[cc]{$x^3_4$}}
\put(25,7){\makebox(0,0)[cc]{$x^1_3$}}
\qbezier(15,10)(38,-5)(60,10)
\end{picture}
\caption{The graph $G$ for the formula $f=C_1\wedge C_2\wedge C_3$ with 
$C_1=x_1\vee x_2\vee x_3$,
$C_2=\bar{x}_1\vee x_4\vee \bar{x}_2$, and 
$C_3=x_1\vee x_3\vee x_4$.}\label{fig3}
\end{center}
\end{figure}

\noindent For the hardness of deciding equality in (\ref{e3}) or (\ref{e3b}),
we describe the efficient construction of a graph $G$ such that 
\begin{eqnarray}\label{ee1}
\mbox{$f$ is satisfiable} 
\,\,\,\,\,\,\,\,\Leftrightarrow \,\,\,\,\,\,\,\,
{\rm diss}(G)=2\alpha(G)
\,\,\,\,\,\,\,\,\Leftrightarrow \,\,\,\,\,\,\,\,
{\rm diss}(G)=2\nu_s(G).
\end{eqnarray}
For every clause $C_i=x\vee y\vee z$ in $f$, 
where $x$, $y$, and $z$ are the three literals in $C_i$, 
we introduce the four vertices $x^i$, $y^i$, $z^i$, and $c^i$ in $G$ 
that induce a clique $G_i$,
where $x^i$, $y^i$, and $z^i$ 
are associated with the three literals $x$, $y$, and $z$ in $C_i$.
Note that $G$ has order $4m$.
For every two vertices $u$ and $v$ belonging to different cliques $G_i$
such that the literal associated with $u$ 
is the negation of the literal associated with $v$,
we add to $G$ the edge $uv$.
This completes the construction of $G$; 
see Figure \ref{fig3} for an illustration.

Clearly, the set $I=\{ c^1,\ldots,c^m\}$ is a maximum independent set of $G$,
in particular, we have $\alpha(G)=m$.
The structure of $G$ easily implies that $G$ 
has a maximum induced matching $M$ 
that only contains edges from $G_1\cup \ldots \cup G_m$,
in fact, any edge in $M$ between a vertex $x$ in $G_i$ and some $G_j$ for $i\not=j$
can be replaced by the edge $xc^i$.
Similarly, the graph $G$ has a maximum dissociation set $D$ 
such that all edges induced by $D$ belong to $G_1\cup \ldots \cup G_m$.
These observations easily imply that $G$ satisfies (\ref{e3d}) with equality,
that is, we have ${\rm diss}(G)=\alpha(G)+\nu_s(G)$.

As observed by Karp, the formula $f$ is satisfiable if and only if 
$G-I$ has an independent set $I'$ of order $m$.
If $f$ is satisfiable, and $I'$ is as above, then 
$I\cup I'$ is a maximum dissociation set in $G$,
and the edges spanned by $I\cup I'$ form a maximum induced matching in $G$,
that is, we have
${\rm diss}(G)=2\nu_s(G)=2m=2\alpha(G)$.
Conversely, if ${\rm diss}(G)=2\alpha(G)$,
then $G$ has a maximum dissociation set $D$ containing $I$,
and $D\setminus I$ 
is an independent set in $G-I$ of order $m$,
that is, it follows that $f$ is satisfiable.
Similarly, if ${\rm diss}(G)=2\nu_s(G)$,
then (\ref{e3d}) implies $\nu_s(G)=m$, and
$G$ has a maximum induced matching $M$ covering $I$,
and the vertices covered by $M$ not in $I$
form an independent set in $G-I$ of order $m$,
that is, again it follows that $f$ is satisfiable.
This completes the proof of (\ref{ee1}),
which shows the NP-hardness of deciding equality in (\ref{e3}) or (\ref{e3b}).

\medskip

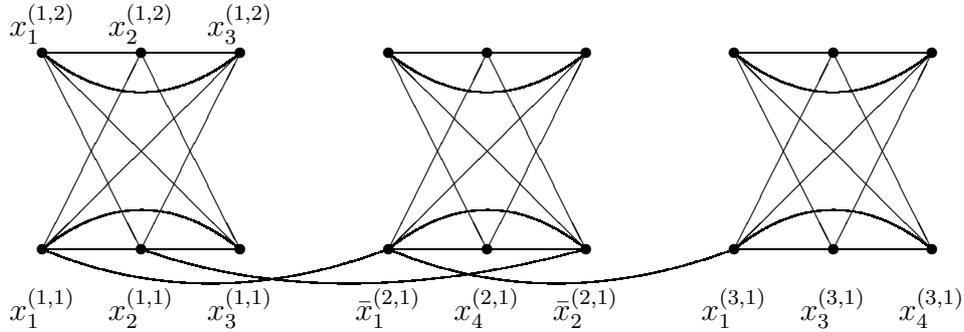
\begin{figure}[h]
\begin{center}
\unitlength 1.3mm 
\linethickness{0.4pt}
\ifx\plotpoint\undefined\newsavebox{\plotpoint}\fi 
\begin{picture}(96,38)(0,0)
\put(5,10){\circle*{1}}
\put(5,30){\circle*{1}}
\put(40,10){\circle*{1}}
\put(40,30){\circle*{1}}
\put(75,10){\circle*{1}}
\put(75,30){\circle*{1}}
\put(15,10){\circle*{1}}
\put(15,30){\circle*{1}}
\put(50,10){\circle*{1}}
\put(50,30){\circle*{1}}
\put(85,10){\circle*{1}}
\put(85,30){\circle*{1}}
\put(25,10){\circle*{1}}
\put(25,30){\circle*{1}}
\put(60,10){\circle*{1}}
\put(60,30){\circle*{1}}
\put(95,10){\circle*{1}}
\put(95,30){\circle*{1}}
\put(5,10){\line(1,0){20}}
\put(5,30){\line(1,0){20}}
\put(40,10){\line(1,0){20}}
\put(40,30){\line(1,0){20}}
\put(75,10){\line(1,0){20}}
\put(75,30){\line(1,0){20}}
\qbezier(25,10)(15,18)(5,10)
\qbezier(25,30)(15,22)(5,30)
\qbezier(60,10)(50,18)(40,10)
\qbezier(60,30)(50,22)(40,30)
\qbezier(95,10)(85,18)(75,10)
\qbezier(95,30)(85,22)(75,30)
\qbezier(5,10)(21.5,3)(40,10)
\qbezier(40,10)(56.5,3)(75,10)
\qbezier(15,10)(34.5,3)(60,10)
\put(5,4){\makebox(0,0)[cc]{$x^{(1,1)}_1$}}
\put(40,4){\makebox(0,0)[cc]{$\bar{x}^{(2,1)}_1$}}
\put(75,4){\makebox(0,0)[cc]{$x^{(3,1)}_1$}}
\put(15,4){\makebox(0,0)[cc]{$x^{(1,1)}_2$}}
\put(60,4){\makebox(0,0)[cc]{$\bar{x}^{(2,1)}_2$}}
\put(50,4){\makebox(0,0)[cc]{$x^{(2,1)}_4$}}
\put(85,4){\makebox(0,0)[cc]{$x^{(3,1)}_3$}}
\put(95,4){\makebox(0,0)[cc]{$x^{(3,1)}_4$}}
\put(25,4){\makebox(0,0)[cc]{$x^{(1,1)}_3$}}
\put(5,33){\makebox(0,0)[cc]{$x^{(1,2)}_1$}}
\put(15,33){\makebox(0,0)[cc]{$x^{(1,2)}_2$}}
\put(25,33){\makebox(0,0)[cc]{$x^{(1,2)}_3$}}
\put(5,30){\line(1,-2){10}}
\put(40,30){\line(1,-2){10}}
\put(75,30){\line(1,-2){10}}
\put(15,10){\line(1,2){10}}
\put(50,10){\line(1,2){10}}
\put(85,10){\line(1,2){10}}
\put(25,30){\line(-1,-1){20}}
\put(60,30){\line(-1,-1){20}}
\put(95,30){\line(-1,-1){20}}
\put(5,10){\line(1,2){10}}
\put(40,10){\line(1,2){10}}
\put(75,10){\line(1,2){10}}
\put(15,30){\line(1,-2){10}}
\put(50,30){\line(1,-2){10}}
\put(85,30){\line(1,-2){10}}
\put(25,10){\line(-1,1){20}}
\put(60,10){\line(-1,1){20}}
\put(95,10){\line(-1,1){20}}
\end{picture}
\caption{The graph $H$ for the formula $f=C_1\wedge C_2\wedge C_3$ with 
$C_1=x_1\vee x_2\vee x_3$,
$C_2=\bar{x}_1\vee x_4\vee \bar{x}_2$, and 
$C_3=x_1\vee x_3\vee x_4$.}\label{fig4}
\end{center}
\end{figure}

\noindent For the hardness of deciding equality in (\ref{e3c}),
we describe the efficient construction of a graph $H$ 
such that $f$ is satisfiable if and only if ${\rm diss}(H)=\alpha(H)$.
For every clause $C_i=x\vee y\vee z$ in $f$, 
where $x$, $y$, and $z$ are the three literals in $C_i$, 
we introduce the six vertices 
$x^{(i,1)}$, $y^{(i,1)}$, $z^{(i,1)}$, $x^{(i,2)}$, $y^{(i,2)}$, and $z^{(i,2)}$ in $H$
that induce a subgraph $H_i$ that is a clique 
minus the three edges $x^{(i,1)}x^{(i,2)}$, $y^{(i,1)}y^{(i,2)}$, and $z^{(i,1)}z^{(i,2)}$.
Similarly as above, the vertices $x^{(i,1)}$, $y^{(i,1)}$, and $z^{(i,1)}$ in $H_i$
are associated with the three literals $x$, $y$, and $z$ in $C_i$.
Note that $H$ has order $6m$.
For every two vertices $u$ and $v$ belonging to different subgraphs $H_i$
that are associated with literals 
such that the literal associated with $u$ 
is the negation of the literal associated with $v$,
we add to $H$ the edge $uv$.
This completes the construction of $H$; 
see Figure \ref{fig4} for an illustration.

Since every dissociation set in $H$ intersects each $H_i$ in at most two vertices,
and selecting two vertices with exponent $(i,2)$ in $H_i$ for each $i$ in $[m]$ 
yields a dissociation set in $H$, we have ${\rm diss}(H)=2m$.
By the structure of $H$, we have ${\rm diss}(H)=\alpha(H)$
if and only if $H$ has an independent set that contains, for every $i$ in $[m]$,
exactly one of the vertices with exponent $(i,1)$.
As noted above, this is equivalent to the satisfiability of $f$,
which shows the NP-hardness of deciding equality in (\ref{e3c}).

\medskip

\noindent For the hardness of deciding equality in (\ref{e3d}),
we describe an efficient reduction from the NP-complete {\sc Independent Set} problem.
Therefore, let $(G,k)$ be an instance of this problem,
that is, the problem of deciding whether $\alpha(G)\geq k$.
Possibly by adding isolated vertices to $G$ 
and increasing $k$ for each added vertex by one, 
we may assume that $2(k-1)>n\geq 2$, 
where $n$ is the order of $G$.
We describe the efficient construction of a graph $H$ 
such that $\alpha(G)\geq k$
if and only if (\ref{e3d}) does {\it not} hold with equality.
The graph $H$ arises from $G$
\begin{itemize}
\item by adding, for every vertex $u$ of $G$, a new vertex $u'$ as well as the edge $uu'$, and 
\item by adding a disjoint copy of $(k-1)K_2$, that is, $k-1$ further independent edges,
as well as all possible edges between the original vertices of $G$ 
and the vertices of the copy of $(k-1)K_2$.
\end{itemize}
If $V$ denotes the vertex set of $G$, then the vertex set of $H$
is $V\cup V'\cup W$, where 
$V'=\{ u':u\in V\}$, 
$W$ is the set of the $2(k-1)$ vertices of the copy of $(k-1)K_2$,
there are all possible edges between $W$ and $V$,
and no edges between $V'$ and $W$.
The order of $H$ is $2n+2(k-1)$.
It is easy to see that $\alpha(H)=n+k-1$,
in fact, the set $V'$ together with one vertex on each of the $k-1$ edges within $W$ 
yields a maximum independent set in $H$.

Our next goal is to show ${\rm diss}(H)=n+2(k-1)$.
Since $V'\cup W$ is a dissociation set in $H$, we have ${\rm diss}(H)\geq n+2(k-1)$.
Now, let $D$ be a maximum dissociation set in $H$.
If $D$ intersects both $V$ and $W$, then $D$ contains exactly 
one vertex from $V$, 
one vertex from $W$, and
all but one vertices from $V'$,
that is, $|D|\leq 1+1+n-1=n+1<n+2(k-1)$, which is a contradiction.
If $D$ does not intersect $W$, then $|D|\leq |V\cup V'|=2n<n+2(k-1)$, 
which is a contradiction.
Thus, the set $D$ does not intersect $V$, which,
by the choice of $D$, implies $D=V'\cup W$,
and, hence, we obtain ${\rm diss}(H)=|D|=|V'\cup W|=n+2(k-1)$ as desired.

In view of the $k-1$ independent edges in $W$, we have $\nu_s(H)\geq k-1$.
Since ${\rm diss}(H)=\alpha(H)+k-1$, 
in order to complete the proof, it suffices to show that
$\alpha(G)\geq k$ if and only if $\nu_s(H)\geq k$.
If $I$ is an independent set in $G$ of order at least $k$,
then $\{ uu':u\in I\}$ is an induced matching in $H$,
hence $\alpha(G)\geq k$ implies $\nu_s(H)\geq k$.
Now, suppose that $\nu_s(H)\geq k$,
and let $M$ be a maximum induced matching in $H$
containing as few edges with both endpoints in $V$ as possible.
If $M$ contains an edge with both endpoints in $W$,
then all edges in $M$ have both endpoints in $W$,
which implies the contradiction $|M|\leq k-1$.
If $M$ contains an edge between $W$ and $V$,
then we obtain the contradiction $|M|=1$.
Hence, no edge in $M$ covers any vertex of $W$.
If $uv\in M$ for $u,v\in V$,
then $M\setminus \{ uv\}\cup \{ uu'\}$ 
is a maximum induced matching in $H$
containing fewer edges with both endpoints in $V$ than $M$.
Hence, the choice of $M$ implies that the set of $|M|\geq k$ vertices from $V$ covered by an edge from $M$ is an independent set in $G$,
that is, $\alpha(G)\geq k$.

This completes the proof of Theorem \ref{theorem2}.

\section{Conclusion}

Our initial motivation to consider the extremal graphs for (\ref{e2}) 
was the attempt to improve the approximation algorithm of Hosseinian and Butenko \cite{hobu}.
This remains to be done.
As explained in the introduction, 
the complexity of recognizing the bipartite extremal graphs  
for (\ref{e3b}), (\ref{e3c}), and (\ref{e3d}) remains open.
We believe that all three problems are hard.
The estimates (\ref{e1}), (\ref{e2}), and (\ref{e3}) imply 
$\alpha(G-M)\leq {\rm diss}(G)\leq 2\alpha(G-M)$
for a given graph $G$ and a given maximum matching $M$ in $G$.
Theorem \ref{theorem2} easily implies that the extremal graphs 
for these two inequalities are also hard to recognize.
In fact, if $G$ is a graph with $\alpha(G)\geq 3$,
the graph $H$ arises from the disjoint union of $G$ and $K_{n(G)}$ 
by adding all possible edges between $V(G)$ and $V(K_{n(G)})$,
and $M$ is a perfect matching of $H$ using only edges between $V(G)$ and $V(K_{n(G)})$,
then 
$\alpha(H)=\alpha(H-M)=\alpha(G)$ and
${\rm diss}(H)={\rm diss}(H-M)={\rm diss}(G)$.
This implies that, for every $c\in \{ 1,2\}$,
the graph $G$ satisfies ${\rm diss}(G)=c\cdot \alpha(G)$
if and only if 
the graph $H$ satisfies 
${\rm diss}(H)=c\cdot \alpha(H-M)$.

\end{document}